\documentclass[preprint,11pt]{elsarticle}

\usepackage[utf8]{inputenc}
\usepackage[T1]{fontenc}

\usepackage{amsmath,amssymb,bm}
\usepackage{graphicx}
\usepackage{float}
\usepackage{booktabs}
\usepackage{booktabs,tabularx,threeparttable}
\usepackage{rotating} 

\usepackage{xcolor}
\definecolor{linkblue}{RGB}{0,95,184}    
\definecolor{urlviolet}{RGB}{140,60,150} 

\usepackage{hyperref}
\hypersetup{
  colorlinks=true,
  linkcolor=linkblue,   
  citecolor=linkblue,   
  urlcolor=urlviolet,   
}

\usepackage{cleveref} 

\journal{.} 
\usepackage{tikz}
\usepackage{pgfplots}
\pgfplotsset{compat=1.18}
\usepackage{siunitx} 
\usepackage{tikz}
\usetikzlibrary{positioning,arrows.meta,decorations.pathreplacing,calc}
\usepackage{tikz}
\usetikzlibrary{shapes.geometric} 
\usepackage{tikz}
\usetikzlibrary{positioning,arrows.meta,fit,backgrounds,decorations.pathreplacing}
\usepackage[a4paper,top=2cm,bottom=2cm,left=1.8cm,right=1.8cm]{geometry}
\usepackage{xcolor}
\definecolor{ElsevierBlue}{RGB}{0,95,184}   
\definecolor{Slate}{RGB}{38,50,56}          
\definecolor{Teal}{RGB}{0,125,125}          
\definecolor{Sky}{RGB}{227,239,252}         
\definecolor{Aqua}{RGB}{204,232,246}        
\definecolor{Indigo}{RGB}{76,81,191}        
\definecolor{Amber}{RGB}{217,119,6}         
\definecolor{Warm}{RGB}{255,247,237}        

\usepackage{tikz}
\usetikzlibrary{
  arrows.meta,
  positioning,
  calc,
  decorations.pathreplacing,
  backgrounds,
  shadows.blur
}

 \usepackage{ltablex}      
 \usepackage{booktabs}     
\usepackage{tabularray}
\UseTblrLibrary{booktabs}
\usepackage{pifont}
\newcommand{\ok}{\ding{51}} 
\newcommand{\x}{\ding{55}}  

 \newcolumntype{L}[1]{>{\raggedright\arraybackslash}p{#1}}
 \newcolumntype{C}[1]{>{\centering\arraybackslash}p{#1}}
 \renewcommand{\arraystretch}{1.12}
 \usepackage[mathlines]{lineno} 

\begin{document}

\begin{frontmatter}

\title{A Selective Review of Modern Stochastic Modeling: SDE/SPDE Numerics, Data-Driven Identification, and Generative Methods with Applications in Biomathematics}

\author[aff1]{Yassine Sabbar}
\author[aff2]{Kottakkaran Sooppy Nisar}

\address[aff1]{IMIA Laboratory, T-IDMS, Department of Mathematics, FST Errachidia, Moulay Ismail University of Meknes, P.O. Box 509, 52000, Boutalamine, Morocco. Emai: \url{y.sabbar@umi.ac.ma}}
\address[aff2]{Department of Mathematics, College of Science and Humanities in Al-Kharj, Prince Sattam bin Abdulaziz University, Al-Kharj 11942, Saudi Arabia. Emai: \url{n.sooppy@psau.edu.sa}}

\begin{abstract}

This review maps 2020-2025 developments in stochastic modeling, highlighting non-standard approaches and their applications to biology and epidemiology. It brings together four strands: (1) core models for systems that evolve with randomness; (2) learning key parts of those models directly from data; (3) methods that can generate realistic synthetic data in continuous time; and (4) numerical techniques that keep simulations stable, accurate, and faithful over long runs. The objective is practical: help researchers quickly see what is new, how the pieces fit together, and where important gaps remain. We summarize tools for estimating changing infection or reaction rates under noisy and incomplete observations, modeling spatial spread, accounting for sudden jumps and heavy tails, and reporting uncertainty in a way that is useful for decisions. We also highlight open problems that deserve near-term attention: separating true dynamics from noise when data are irregular; learning spatial dynamics under random influences with guarantees of stability; aligning training with the numerical method used in applications; preserving positivity and conservation in all simulations; reducing cost while controlling error for large studies; estimating rare but important events; and adopting clear, comparable reporting standards. By organizing the field around these aims, the review offers a concise guide to current methods, their practical use, and the most promising directions for future work in biology and epidemiology.

\end{abstract}


\end{frontmatter}

\section{Introduction}
\label{sec:intro}

Stochastic modeling offers a mathematically rigorous and physically interpretable framework for analyzing dynamical systems whose evolution is influenced by intrinsic randomness \cite{0a}, unresolved multiscale interactions \cite{0b}, or incomplete specification of the governing mechanisms~\cite{1}. In contrast to deterministic formulations (where the system trajectory is uniquely determined by a prescribed set of ordinary or partial differential equations \cite{0c}) stochastic models incorporate random variables or stochastic processes directly into the governing equations~\cite{2}. This probabilistic formulation enables explicit quantification of epistemic and aleatory uncertainty, reproduces variability across independent realizations \cite{2a}, and captures non-deterministic phenomena, such as intermittency, heavy-tailed fluctuations, or regime switching, that lie beyond the representational capacity of purely deterministic models~\cite{3}.

For finite-dimensional systems, the canonical representation is the stochastic differential equation (SDE):
\begin{equation}
    \mathrm{d}X_t = f(X_t,t)\,\mathrm{d}t + G(X_t,t)\,\mathrm{d}W_t,
    \label{eq:SDE}
\end{equation}
where \(X_t \in \mathbb{R}^d\) denotes the state vector, \(f : \mathbb{R}^d \times [0,\infty) \to \mathbb{R}^d\) is the drift function describing the deterministic dynamics, \(G : \mathbb{R}^d \times [0,\infty) \to \mathbb{R}^{d\times m}\) is the diffusion coefficient matrix, and \(W_t\) is an \(m\)-dimensional standard Wiener process modeling temporally continuous Gaussian perturbations~\cite{4}. The drift term \(f\) typically encodes system-specific physical laws or conservation principles, whereas the diffusion term \(G\) parameterizes the intensity, anisotropy, and possible state-dependence of the noise \cite{4a}.  When the state space is infinite-dimensional, one arrives at stochastic partial differential equations (SPDEs), where the solution \(X(t,\xi)\) depends on both temporal (\(t\)) and spatial (\(\xi\)) variables, and the stochastic forcing may act locally or globally in space~\cite{5}. SPDEs arise naturally in the modeling of spatially distributed systems subject to random forcing, such as turbulent flows \cite{5a}, population dispersal \cite{5b}, or reaction–diffusion processes with environmental fluctuations \cite{5c}.

The classical Itô calculus framework has been generalized to accommodate diverse noise structures~\cite{6}, including:  
(i) \emph{Jump–diffusion processes}~\cite{7}, where compensated Poisson random measures introduce finite-activity discontinuities representing abrupt or rare events~\cite{8};  
(ii) \emph{Lévy-driven systems}~\cite{9}, which allow for heavy-tailed increments and infinite-activity jump processes~\cite{10}, capturing burst-like dynamics and extreme variability ~\cite{11}; and  
(iii) \emph{Fractional Brownian motion} and \emph{Volterra-type kernels}~\cite{11a,12}, which encode long-range temporal dependence, self-similarity, and memory effects \cite{12a}.  
In many applications, the diffusion term \(G(X_t,t)\) is itself state-dependent, creating a multiplicative noise structure and introducing bidirectional coupling between the system dynamics and the stochastic perturbations~\cite{13}.

From a statistical perspective, stochastic models serve both as \emph{generative mechanisms} and as \emph{inferential frameworks}~\cite{14,15}. As generative models, they produce synthetic sample paths that preserve empirical distributional properties \cite{15a}, temporal correlation structures \cite{15b}, and extreme-event statistics observed in experimental or observational datasets~\cite{16}. As inferential frameworks, they provide a principled foundation for recovering latent parameters, hidden states \cite{16a}, and entire posterior distributions from noisy, incomplete, or irregularly sampled measurements~\cite{17}.

Parameter estimation for SDEs and SPDEs is supported by a broad methodological arsenal, including maximum likelihood estimation~\cite{17a}, the generalized method of moments~\cite{17b}, and Bayesian inference frameworks~\cite{17c} that yield full posterior distributions and enable comprehensive uncertainty quantification \cite{18}. Computational advances now permit the application of Markov chain Monte Carlo (MCMC) methods (such as Hamiltonian Monte Carlo and particle MCMC~\cite{18a}) to efficiently explore high-dimensional posterior spaces~\cite{19}. For online estimation in systems with streaming data, sequential Monte Carlo methods, such as the particle filter and the ensemble Kalman filter, have proven indispensable~\cite{20}. In scenarios where the likelihood function is analytically intractable, likelihood-free inference techniques, including approximate Bayesian computation and synthetic likelihood methods \cite{20a}, offer practical alternatives by relying on forward simulation and carefully chosen summary statistics in place of explicit likelihood evaluation~\cite{21}.

Recent theoretical advances in high-frequency sampling asymptotics have yielded consistent and asymptotically efficient estimators for both drift and diffusion coefficients under substantially weaker smoothness and regularity conditions than those previously required \cite{22}. These results hold even under model misspecification \cite{23}, endogenous observation noise \cite{24}, or irregular sampling schedules, thereby extending their applicability to empirical datasets where idealized sampling assumptions are violated \cite{25}. In the context of SPDEs, spectral estimation methods (combined with finite-dimensional Galerkin projections \cite{26}) enable accurate recovery of dynamical modes from partial, noisy, and spatially sparse measurements \cite{27}. Parallel developments in computational statistics have given rise to hybrid methodologies that integrate stochastic modeling with machine learning  \cite{28}. In particular, Neural SDE frameworks  \cite{29}, physics-informed neural operators (PINO) \cite{30}, and operator-learning architectures can now recover the functional structure of the drift \(f\) and diffusion \(G\) in~\eqref{eq:SDE} directly from data, while rigorously enforcing physical and statistical constraints \cite{31}. This fusion of stochastic analysis, numerical approximation, and data-driven inference represents a decisive step toward predictive modeling of complex, high-dimensional, and uncertainty-dominated systems \cite{32}.

Over the past decade, several research frontiers have expanded both the theoretical foundations and the computational reach of stochastic modeling \cite{32a,33}. One prominent direction addresses the treatment of high-dimensional and multiscale systems \cite{34}, where the number of interacting stochastic degrees of freedom can be in the hundreds or thousands \cite{35}, as encountered in turbulent fluid dynamics \cite{36}, molecular simulations \cite{37}, and eddy-resolving climate models \cite{38}. In such regimes, direct time integration of the full stochastic system is often computationally prohibitive \cite{39}. To overcome this barrier, dimensionality-reduction techniques (such as dynamic mode decomposition, stochastic principal component analysis, and manifold learning \cite{40}) combined with reduced-order and homogenized models \cite{41}, have enabled the derivation of effective lower-dimensional stochastic dynamics that preserve essential statistical, spectral, and dynamical properties of the original system while significantly reducing computational cost \cite{42}.

Another active research area concerns non-Lipschitz and stiff stochastic dynamics \cite{43}, which arise in systems with super-linear drifts, polynomial nonlinearities, or singular coefficients \cite{44}. Classical explicit discretizations, including the Euler–Maruyama scheme \cite{45}, may fail in such settings due to numerical instability or divergence. To address these challenges, tamed and balanced numerical integrators have been developed to control the growth of nonlinear terms without sacrificing strong or weak convergence properties \cite{46}. In addition, implicit–explicit (IMEX) schemes \cite{47}, exponential integrators \cite{48}, and stochastic Rosenbrock-type methods \cite{49} have been successfully adapted to handle severe stiffness and disparate temporal scales, ensuring stability while maintaining high-order accuracy \cite{50}.

The modeling of heavy-tailed fluctuations and discontinuous trajectories has likewise progressed significantly \cite{51}. Notable developments include SDEs driven by $\alpha$-stable Lévy processes \cite{52}, which generate infinite-variance increments and capture heavy-tailed statistics \cite{53}; compound Poisson processes \cite{54}, which model finite-activity but potentially large-amplitude jumps \cite{55}; and state-dependent jump intensities \cite{56}, allowing the jump rate to vary adaptively with the evolving system state \cite{57}. Such models have become indispensable in quantitative finance \cite{58}, where they describe market shocks and volatility clustering; in reliability engineering \cite{59}, where they model catastrophic system failures; and in epidemiology, where they capture abrupt changes in transmission due to superspreading events or public-health interventions \cite{60}.

Memory-dependent stochastic processes constitute another rapidly advancing frontier \cite{61}. Fractional SDEs and SPDEs \cite{62}, incorporating Caputo or Riemann–Liouville derivatives \cite{63}, accurately model anomalous diffusion \cite{64}, viscoelastic response \cite{65}, and long-range temporal correlations observed in diverse physical, biological \cite{66}, and socioeconomic systems \cite{67}. The inherently nonlocal structure of fractional operators poses severe computational challenges \cite{68}, which have been mitigated by fast convolution algorithms \cite{69}, kernel compression strategies \cite{70}, and spectral approximations \cite{71}, thereby reducing memory and time complexity from \(O(N^2)\) to nearly \(O(N \log N)\) without significant loss of accuracy \cite{72}.

The integration of stochastic modeling with modern machine learning has given rise to a new class of data-driven inference methods \cite{73}. Neural SDEs parameterize the drift and diffusion terms using deep neural networks \cite{74}, enabling direct estimation from raw, noisy, or irregularly sampled data while preserving stochastic well-posedness \cite{75}. Operator-learning paradigms \cite{76}, including physics-informed neural operators, generalize this framework to SPDEs, learning mappings between infinite-dimensional function spaces under physical and statistical constraints \cite{77}. Furthermore, diffusion-based generative models reinterpret high-dimensional data synthesis as the numerical solution of forward–reverse SDEs, providing a probabilistic bridge between stochastic analysis and state-of-the-art generative learning \cite{78,79}.

Advances in statistical inference have paralleled these methodological innovations \cite{80}. Scalable Bayesian inference for SDE and SPDE models now leverages surrogate modeling \cite{81}, Gaussian process emulators \cite{82}, and multi-fidelity Monte Carlo schemes to accelerate likelihood evaluations \cite{83}. This has enabled full posterior sampling in previously intractable high-dimensional settings, thereby facilitating rigorous uncertainty quantification in complex stochastic systems \cite{84}. When exact likelihood evaluation is impossible, likelihood-free approaches such as approximate Bayesian computation (ABC) \cite{85} and synthetic likelihood estimation have emerged as powerful alternatives, relying on forward simulation combined with carefully chosen summary statistics \cite{86}.

The impact of these developments is evident across multiple disciplines. In epidemiology, stochastic compartmental models augmented with demographic noise \cite{87}, Lévy-driven transmission rates \cite{88}, and adaptive behavioral feedback mechanisms have achieved superior realism in reproducing epidemic variability \cite{89}, persistence probabilities, and extinction events \cite{90}. In quantitative finance, regime-switching jump–diffusion models accurately capture sudden volatility spikes and heavy-tailed return distributions \cite{91,92}. In climate science, SPDEs driven by spatially correlated colored noise provide a principled framework for representing structured uncertainty in coupled ocean–atmosphere systems \cite{93}. In robotics and control, stochastic optimal control problems (often formulated via backward SDEs) enable robust decision-making under uncertainty, supporting safe trajectory planning and real-time risk-aware operations in dynamic environments \cite{94}.

This review delivers a synthesis of recent advances in stochastic modeling with a strong emphasis on the dynamic interplay between theoretical developments, statistical inference, and computational algorithm design~\cite{95}. While the mathematical foundations of SDEs and SPDEs are well established, the period 2020--2025 has witnessed a rapid expansion driven by data-driven system identification, the integration of modern machine-learning architectures~\cite{96}, and the development of robust numerical schemes tailored to nonlinearities and stochastic perturbations. These developments are organized into a structured taxonomy that links established stochastic theory with modern approaches shaped by high-dimensional data analysis and physics-informed learning, covering SDE/SPDE numerics, data-driven identification via neural SDEs and filtering, and generative as well as reverse-time methods for stochastic dynamics. The survey critically evaluates numerical integration strategies in terms of stability under stiffness, convergence under non-Lipschitz growth, and fidelity to invariant measures in the presence of jump noise, heavy-tailed perturbations, and fractional effects. It also consolidates recent advances in statistical inference for contemporary stochastic systems, addressing the challenges of estimation and uncertainty quantification in high-dimensional states with partial observations. A distinctive feature of this review is its application-oriented synthesis, mapping the surveyed methods to concrete problems in biology and epidemiology, such as the estimation of time-varying transmission or reaction rates, epidemic forecasting under reporting biases, modeling spatial spread through SPDE formulations, and quantifying uncertainty for decision-making. By compiling and analyzing recent works from 2020 to 2025, this review aims to provide researchers with a clear view of the current state of the art, highlight methodological novelties, and identify promising directions for advancing theoretically rigorous, computationally scalable, and practically impactful models in biological and epidemiological systems.

The remainder of this article is structured to provide a coherent progression from foundational concepts to advanced applications. Section~\ref{sec:learning} examines data-driven identification strategies for stochastic systems, with particular attention to neural SDE frameworks and operator-learning methodologies that enable the direct inference of drift and diffusion structures from data. Section~\ref{sec:generative} explores the emerging role of SDE-based formulations in generative modeling, highlighting their connection to forward–reverse processes and diffusion-based learning. Section~\ref{sec:numerics} offers a detailed survey of numerical integration techniques for SDEs and SPDEs, emphasizing stability, convergence, and long-time accuracy in challenging regimes.  Representative applications in biology and epidemiology  are presented in Section~\ref{sec:applications}, and the review concludes in Section~\ref{sec:outlook} with a synthesis of open research challenges and promising future directions.

\section{Learning Stochastic Dynamics from Data}
\label{sec:learning}

The growing availability of high-resolution, large-scale datasets has shifted the emphasis in stochastic modeling from prescribing dynamics solely on the basis of first-principles derivations to inferring them directly from empirical observations. Within this data-centric paradigm, two complementary methodologies have emerged. The first addresses finite-dimensional systems through neural SDEs \cite{97}, in which both drift and diffusion structures are represented by flexible, trainable mappings estimated from data. The second focuses on infinite-dimensional systems via operator-learning approaches for SPDEs \cite{98}, which aim to approximate stochastic evolution operators in a mesh-independent fashion. Both frameworks are designed to assimilate sparse or noisy measurements, generalize robustly across dynamical regimes, and preserve the analytical properties imposed by the underlying stochastic calculus.

\subsection{Neural SDEs and Drift–Diffusion Identification}
\label{subsec:neural_sde}

In the neural SDE setting, the drift and diffusion terms are replaced by parameterized functions learned from data, while the governing Itô form is retained \cite{74}. A general formulation on $\mathbb{R}^d$ is
\begin{equation}
\mathrm{d}X_t = f_\theta(X_t,t)\,\mathrm{d}t + G_\theta(X_t,t)\,\mathrm{d}W_t,
\qquad X_0 \sim \mu_0,
\label{eq:neuralSDE}
\end{equation}
where $X_t\in\mathbb{R}^d$ is the state vector at time $t$, $f_\theta:\mathbb{R}^d\times[0,T]\to\mathbb{R}^d$ is the drift, $G_\theta:\mathbb{R}^d\times[0,T]\to\mathbb{R}^{d\times m}$ is the diffusion coefficient matrix, $\theta$ denotes all trainable parameters, $W_t$ is an $m$-dimensional standard Wiener process, and $\mu_0$ is the distribution of the initial state. When discrete samples $\{x_{t_k}\}_{k=0}^n$ with time steps $\Delta t_k=t_{k+1}-t_k$ are available, one common identification approach is maximum likelihood under a discretized transition model \cite{100}. Using the Euler–Maruyama approximation \cite{101}, the conditional distribution is
\[
X_{t_{k+1}}\mid X_{t_k} \approx \mathcal{N}\!\Big(x_{t_k} + f_\theta(x_{t_k},t_k)\,\Delta t_k,\; \Sigma_\theta(x_{t_k},t_k)\,\Delta t_k\Big),
\]
where $\Sigma_\theta := G_\theta G_\theta^\top$ is the diffusion covariance. Maximization of this Gaussian likelihood yields an estimate of $\theta$ \cite{102}. For state-dependent diffusion, higher-order schemes such as Milstein or Itô–Taylor expansions improve accuracy \cite{103}. Positive semidefiniteness of $\Sigma_\theta$ is enforced via Cholesky factorization, and drift growth is regulated through tamed or monotonicity-preserving architectures to ensure well-posedness and strong solution uniqueness \cite{104}.

An alternative approach trains in continuous time by differentiating through the SDE solver. Fixing noise seeds $\xi_k\sim\mathcal{N}(0,I_m)$, a single solver step is denoted by
\[
X_{t_{k+1}} = \Phi_\theta(X_{t_k},t_k;\,\xi_k),
\]
where $\Phi_\theta$ encodes the numerical update rule \cite{105}. Gradients $\nabla_\theta \mathbb{E}[\mathcal{L}(X_{0:n})]$ are obtained via direct backpropagation or through the adjoint SDE, the latter reducing memory from $\mathcal{O}(n)$ to $\mathcal{O}(1)$ in the number of steps \cite{106}. For partially observed systems $Y_k = H X_{t_k} + \varepsilon_k$, where $H$ is the observation matrix and $\varepsilon_k$ is measurement noise, the neural SDE can be coupled with a differentiable filter (such as an ensemble Kalman update or amortized particle weighting) and trained by maximizing a filtering log-likelihood or a variational lower bound \cite{107}.

When likelihood-based methods are unreliable, for example under irregular sampling or unknown observation noise, moment-based estimation offers a robust alternative \cite{108}. For any smooth function $g:\mathbb{R}^d\to\mathbb{R}$, Itô’s formula gives the local martingale:
\[
M_t^g = g(X_t) - g(X_0) - \int_0^t \Big(\nabla g^\top f_\theta + \tfrac12 \operatorname{Tr}[\Sigma_\theta \nabla^2 g]\Big)(X_s,s)\,\mathrm{d}s.
\]
Enforcing that empirical averages of $M_{t_k}^g$ vanish yields unbiased estimating equations even with noisy or irregularly spaced data \cite{109}. The diffusion can also be estimated from the quadratic variation:
\[
\frac{1}{\Delta t_k}(X_{t_{k+1}}-X_{t_k})(X_{t_{k+1}}-X_{t_k})^\top \approx \Sigma_\theta(X_{t_k},t_k),
\]
with bias corrections applied for finite $\Delta t_k$ \cite{110}. In models with jumps, drift inference may employ Girsanov’s theorem \cite{111}, while the jump component can be estimated from its Lévy–Khintchine representation \cite{112}, using, for example, compound Poisson thinning or variational approximations for infinite-activity cases \cite{113}.

Careful treatment of identifiability is essential: without constraints, drift effects can be spuriously absorbed into the diffusion term \cite{114}. Common remedies include imposing structural restrictions on $G_\theta$ (e.g., diagonal or low-rank form), penalizing deviations $\|\widehat{\Sigma} - \Sigma_\theta\|$ over short time increments, or leveraging stationary distributions \cite{115}. If $\pi$ is the invariant measure of \eqref{eq:neuralSDE}, $f_\theta$ and $\Sigma_\theta$ satisfy the stationary Fokker–Planck equation
\[
\nabla\!\cdot\!\big(f_\theta \pi\big) - \tfrac12 \sum_{i,j}\partial_{i}\partial_{j}\!\big(\Sigma_{\theta,ij}\pi\big) = 0,
\]
which can serve as a weak constraint when long-run data are available \cite{116}. In practice, training objectives often combine the transition likelihood, quadratic-variation penalties, weak Fokker–Planck constraints, and stability regularizers that enforce linear growth bounds or one-sided Lipschitz conditions to promote ergodicity \cite{117}.

Uncertainty quantification may be addressed via Bayesian neural SDEs \cite{118}, where priors are placed on $\theta$ and posteriors are explored using stochastic-gradient MCMC or ensemble variational inference \cite{119}, or by nonparametric bootstrap resampling over both trajectories and noise seeds \cite{120}. Robustness to out-of-distribution states can be enhanced by spectral normalization in the drift and by imposing lower bounds on diffusion coefficients \cite{121}. In controlled systems with exogenous inputs $u_t$, the drift $f_\theta(x,t,u)$ can be parameterized in a control-affine form, enabling integration with stochastic optimal control through Hamilton–Jacobi–Bellman surrogates or pathwise policy-gradient methods \cite{122}.

\subsection{Operator Learning for SPDEs (PINO, Neural Operators)}
\label{subsec:operator_learning}

Stochastic partial differential equations on spatial domains $\Omega\subset\mathbb{R}^p$ can be expressed in abstract form as
\begin{equation}
\mathrm{d}u(t) = \mathcal{A}_{\phi}\big(u(t)\big)\,\mathrm{d}t + \mathcal{B}_{\phi}\big(u(t)\big)\,\mathrm{d}W_t,
\qquad u(0)=u_0,\quad t\in[0,T],
\label{eq:SPDE}
\end{equation}
where $u(t):\Omega\to\mathbb{R}^q$ denotes the state variable at time $t$, $\mathcal{A}_\phi$ is a (possibly nonlinear) differential operator parameterized by $\phi$, $\mathcal{B}_\phi$ is a noise-multiplication operator, and $W_t$ is a cylindrical Wiener process on an appropriate Hilbert space. The initial condition is $u_0\in\mathcal{H}$, where $\mathcal{H}$ is the state space \cite{123}. The learning objective is to approximate the solution operator
\[
\mathcal{S} : (u_0, f, \xi) \mapsto u(\cdot),
\]
which maps initial data $u_0$, deterministic forcing terms $f:\Omega\times[0,T]\to\mathbb{R}^q$, and stochastic realizations $\xi$ to full trajectories $u$ \cite{124}. Neural operator architectures approximate $\mathcal{S}$ directly in function space, enabling mesh-independent generalization and evaluation across discretizations \cite{77}. Two principal designs have been particularly effective.

The first class, Fourier or kernel neural operators, alternates local linear mappings with global integral transforms of the form
\[
v_{k+1} = \sigma\!\big(W v_k + \mathcal{K}_\psi v_k\big),
\]
where $v_k$ denotes the feature representation at layer $k$, $W$ is a local (pointwise) linear transformation, $\sigma$ is a nonlinear activation, and $\mathcal{K}_\psi$ is a global integral operator with kernel $K_\psi(x,y)$ parameterized either in physical space or via truncated Fourier multipliers \cite{126}. Stacking such layers yields a mapping $\mathcal{G}_\Psi : \mathcal{X} \to \mathcal{Y}$ that approximates $\mathcal{S}$. In the stochastic setting, randomness enters through sampled forcings $f$ or noise realizations $\xi$, and the training objective minimizes trajectory- or distribution-based losses over ensembles of noise samples \cite{127}. Stability is promoted by constraining spectral multipliers (e.g., bounding spectral radius) and incorporating skip connections designed to mimic the dissipative structure of the underlying PDE \cite{128}.

The second class, PINO, augments data-driven learning with weak enforcement of the governing equations in \eqref{eq:SPDE} \cite{129}. Given a set of test functions $\{\varphi_j\}_{j=1}^J\subset\mathcal{H}$, the variational residual is defined by
\[
\mathcal{R}_\phi(u) := \sum_{j=1}^J \Big\langle \partial_t u - \mathcal{A}_\phi(u),\;\varphi_j \Big\rangle,
\]
where $\langle \cdot, \cdot \rangle$ denotes the $L^2(\Omega)$ inner product. Stochastic forcing is incorporated either by sampling noise paths $\xi$ and enforcing Itô weak identities, or by matching statistical quantities (e.g., covariances, energy spectra) derived from theory \cite{130}. The composite training objective takes the form
\[
\min_{\Psi,\phi}\; \mathbb{E}_{\xi}\!\left[\| \mathcal{G}_\Psi(u_0,f,\xi)-u \|_{\mathcal{Y}}^2\right]
+ \alpha\,\mathbb{E}_{\xi}\!\left[\|\mathcal{R}_\phi(\mathcal{G}_\Psi(u_0,f,\xi))\|^2\right]
+ \beta\,\mathcal{J}_{\text{BC/IC}},
\]
where $\|\cdot\|_{\mathcal{Y}}$ denotes the $L^2$ norm over $\Omega\times[0,T]$, $\mathcal{J}_{\text{BC/IC}}$ enforces boundary and initial conditions, and $\alpha,\beta>0$ are weighting parameters \cite{131}. This formulation reduces data requirements and mitigates nonphysical extrapolation in regimes where stochastic forcing excites unresolved spatial or temporal scales.

Accurate noise representation is critical. A common approach is to expand $\xi(t,x)$ in a Karhunen–Loève basis $\{\phi_i(x)\}_{i=1}^r$ \cite{132},
\[
\xi(t,x) = \sum_{i=1}^r a_i(t)\,\phi_i(x),
\]
where $a_i(t)$ are temporally stochastic coefficients and the $\phi_i$ are spatial modes, possibly learned jointly with the operator \cite{133}. For multiplicative noise, maintaining Itô–Stratonovich consistency requires either explicit correction terms in the residual or learned commutator terms that approximate the conversion \cite{134}. If an invariant measure $\pi$ is known or estimable, an auxiliary ergodic loss
\[
\mathbb{E}_{\pi}[\mathcal{Q}(u)] \approx \frac{1}{T}\int_0^T \mathcal{Q}\big(\mathcal{G}_\Psi(u_0,f,\xi)(t)\big)\,\mathrm{d}t
\]
is used to align surrogate long-time statistics with those of the true system, where $\mathcal{Q}$ denotes a chosen diagnostic such as energy, enstrophy, or a spectral functional \cite{135}.

Generalization across meshes and domains is supported by spectral parameter tying to ensure resolution scalability, coordinate-free message passing for unstructured meshes, and discrete conservation constraints such as divergence-free enforcement for incompressible flows or positivity for density fields \cite{136}. For stiff or multiscale SPDEs, multi-rate operator architectures couple coarse global updates with local fine-scale correctors, analogous to heterogeneous multiscale methods \cite{137}. Uncertainty quantification can be incorporated via Bayesian priors on operator weights, ensembles over noise realizations \cite{138}, or linearization-based covariance approximations using Gauss–Newton methods \cite{139}.

A hybrid approach couples finite-dimensional latent SDEs with learned spatial bases, representing the solution as $u(t,x)\approx \sum_{i=1}^r c_i(t)\,\phi_i(x)$, where $c(t)\in\mathbb{R}^r$ evolves under a neural SDE and $\{\phi_i\}$ are trainable spatial modes \cite{140}. This yields interpretable low-rank stochastic dynamics with mesh-independent inference and efficient sampling, applicable to tasks such as stochastic control, data assimilation, and probabilistic forecasting \cite{141}.
\section{Generative Modeling as SDEs}
\label{sec:generative}

\begin{figure}[H]
    \centering
    \includegraphics[width=1\textwidth,keepaspectratio]{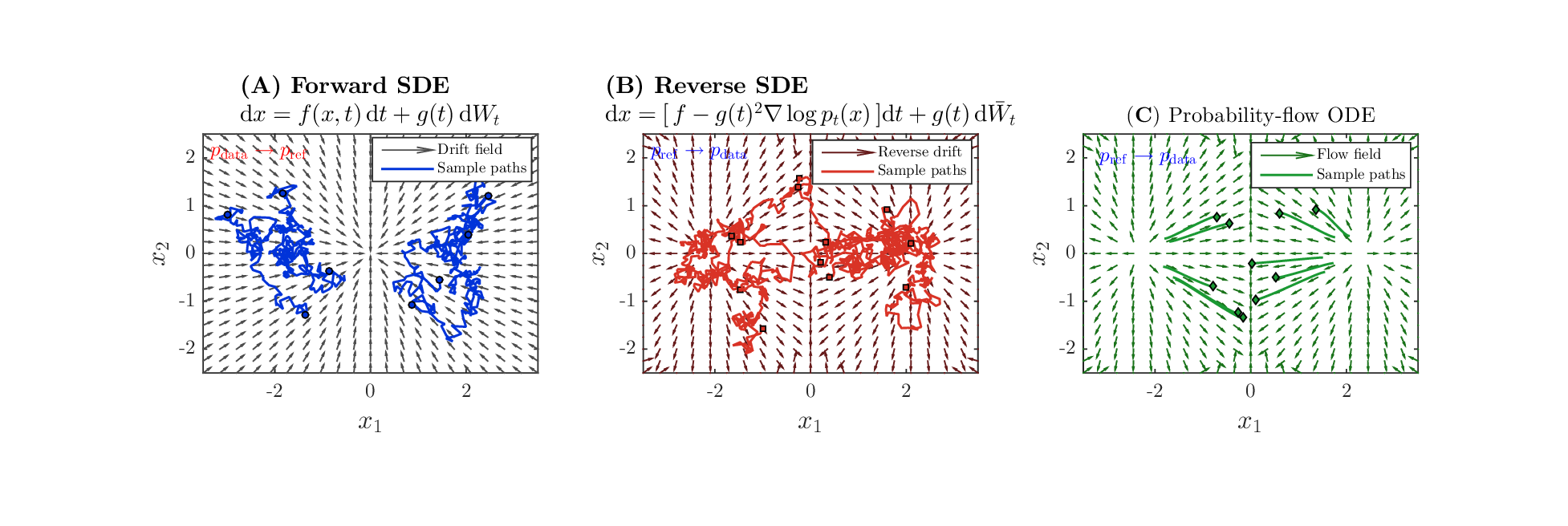}
    \vspace{-1.5cm}
    \caption{
        Forward–reverse SDE framework and its deterministic probability-flow ODE counterpart. 
        Panel~(A) shows the forward SDE transporting $p_{\mathrm{data}}$ to a reference Gaussian $p_{\mathrm{ref}}$. 
        Panel~(B) shows the reverse SDE mapping $p_{\mathrm{ref}}$ back to $p_{\mathrm{data}}$ via the learned score function $\nabla \log p_t(x)$. 
        Panel~(C) depicts the equivalent probability-flow ODE sharing the same marginals as the forward and reverse SDEs but evolving deterministically. 
        Vector fields represent drift components; colored trajectories indicate sample paths. 
            }
    \label{fig:forward_reverse_sde}
\end{figure}

Generative modeling aims to construct stochastic transformations that map a tractable reference distribution onto a complex target distribution \cite{14}. Recent advances have established that score-based generative models and denoising diffusion models can be formulated entirely within the stochastic differential equation (SDE) framework introduced in \eqref{eq:SDE} \cite{78}. In this setting, a \emph{forward} process progressively perturbs samples from the data distribution $p_{\mathrm{data}}$ into a simple reference law $p_{\mathrm{ref}}$, typically a standard Gaussian, over a fixed time horizon  \cite{144}. A \emph{reverse-time} stochastic process then reconstructs new samples by inverting this diffusion \cite{145}, as illustrated in Fig.~\ref{fig:forward_reverse_sde}.

Let $p_t$ denote the marginal distribution of the forward process at time $t$. Under regularity conditions, Anderson’s time-reversal theorem \cite{146} yields the reverse-time SDE:
\begin{equation}
\mathrm{d}x(t) = \Big[ f(x(t),t) - g(t)^{2} \, \nabla_x \log p_t\big(x(t)\big) \Big]\,\mathrm{d}t 
+ g(t)\,\mathrm{d}\bar{w}(t),
\label{eq:reverseSDE}
\end{equation}
where $\bar{w}(t)$ is a standard Wiener process evolving backward in time, $f$ and $g$ are the drift and diffusion coefficients of the forward dynamics, and $\nabla_x \log p_t$ is the score function of $p_t$. Since $p_t$ is not known in closed form, the score function is approximated by a trainable surrogate $\hat{s}_\varphi(x,t)$, estimated through score-matching or denoising-score-matching objectives, ensuring asymptotic consistency in the limit of infinite data and model capacity \cite{147}.

Numerical generation is performed by discretizing \eqref{eq:reverseSDE} with schemes such as Euler–Maruyama \cite{148}, stochastic Runge–Kutta \cite{149}, or predictor–corrector integrators \cite{150}. The choice of scheme directly influences sampling bias and variance, with high-order solvers yielding greater fidelity for a given computational budget \cite{151}. An important alternative is the \emph{probability-flow ordinary differential equation} associated with the forward process, which evolves deterministically yet preserves the same marginal distributions as the reverse SDE, often permitting larger stable time steps \cite{152}.

In scientific applications, the forward dynamics can be tailored to preserve statistical invariants, enforce spectral constraints, or emulate dissipative mechanisms from coarse-grained physics \cite{153}. The learned reverse process then reconstructs fine-scale structure and corrects higher-order moments, enabling physically consistent synthesis in high-dimensional contexts such as turbulent flows, molecular dynamics, and climate modeling \cite{154}. Computational efficiency can be further enhanced through adaptive time stepping, variance-reduction strategies (e.g., control variates), and reduced-order latent representations \cite{155}.

A natural generalization introduces exogenous control inputs $u_t$ into the drift, yielding a controlled reverse-time process of the form 
\[
f(x,t,u_t) - g(t)^{2} \, \nabla_x \log p_t(x,t),
\] 
which connects generative modeling with stochastic optimal control \cite{156}. In this interpretation, the reverse dynamics act as an optimal feedback policy minimizing a cost functional equivalent to the negative log-likelihood, enabling targeted sampling subject to trajectory or terminal constraints \cite{157}. This formulation has direct implications for data assimilation, scenario generation, and constraint-aware simulation \cite{158}.

By embedding generative modeling entirely within the SDE formalism and incorporating domain-specific constraints, one obtains a unifying mathematical structure capable of preserving statistical integrity while extending these methods from static data generation to the synthesis and control of complex stochastic dynamical systems.

\section{Numerical Methods for SDEs and SPDEs}
\label{sec:numerics}

The accurate and efficient numerical integration of stochastic differential equations and their infinite–dimensional counterparts is central to both predictive modeling and statistical inference. The numerical discretization must preserve the fundamental properties of the underlying stochastic dynamics, including convergence order, stability, invariant measures, and qualitative behavior over long time horizons. This section surveys recent methodological developments across four interconnected areas.

\subsection{Strong and Weak Convergence, and Stability}
\label{subsec:convergence_stability}

\begin{figure}[htbp]
    \centering
    \includegraphics[width=0.6\textwidth,keepaspectratio]{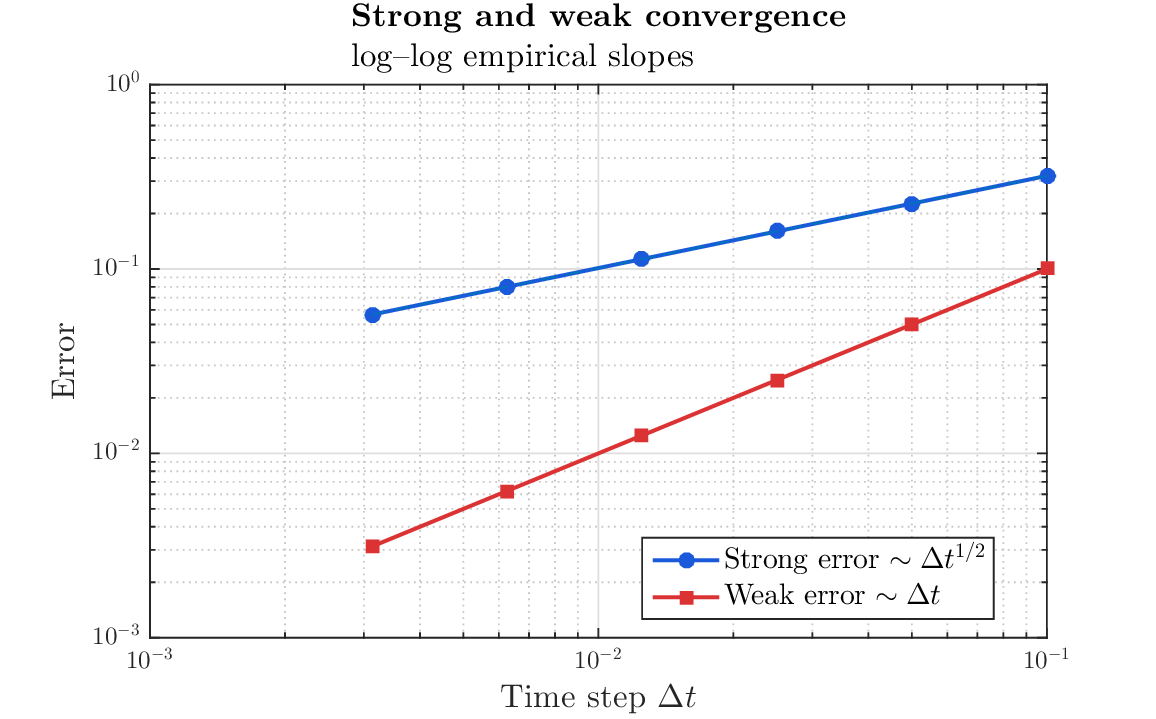}
    \caption{Strong and weak convergence on a log–log scale for representative schemes. The strong error exhibits the expected order one half for an Euler–Maruyama discretization with multiplicative noise, whereas the weak error attains order one for smooth test functionals. Dashed guides indicate the asymptotic slopes used to assess rates. The same methodology extends to semidiscrete SPDEs, where temporal and spatial errors interact through the smoothing properties of the linear part and the regularity of the stochastic forcing.}
    \label{fig:conv_strong_weak}
\end{figure}

In the numerical treatment of stochastic differential equations, strong convergence quantifies the mean–square error between the numerical approximation $X^{\Delta t}_T$ and the exact solution $X_T$ at a fixed terminal time $T$ \cite{159},
\[
\big( \mathbb{E}\,\|X_T - X^{\Delta t}_T\|^2 \big)^{1/2} = \mathcal{O}\big((\Delta t)^\gamma\big),
\]
while weak convergence measures the error in expectations of sufficiently smooth observables $\varphi$,
\[
\big| \mathbb{E}\,\varphi(X_T) - \mathbb{E}\,\varphi(X^{\Delta t}_T) \big| = \mathcal{O}\big((\Delta t)^\beta\big).
\]
For drift–diffusion systems with non–globally Lipschitz drift, polynomial growth, or multiplicative noise, modern analyses derive orders $\gamma$ and $\beta$ under one-sided Lipschitz and coercivity conditions, often with tamed, balanced, or implicit updates to ensure well-posedness of the discrete dynamics \cite{43,160,162}. In infinite dimensions, the same notions apply to semidiscrete SPDEs in Hilbert spaces, where temporal error interacts with spatial discretization error through the smoothing properties of the linear generator and the regularity of the noise \cite{163}. Stability considerations extend beyond mean-square boundedness to include almost-sure exponential stability  \cite{164}, boundedness in probability, and preservation of qualitative invariants over long horizons  \cite{165}. Empirical slopes for representative strong and weak errors, and their asymptotic rates, are illustrated in Fig.~\ref{fig:conv_strong_weak}.

\subsection{Tamed, Implicit, and Balanced Schemes}
\label{subsec:tamed_implicit}

Classical explicit schemes may diverge when applied to SDEs with superlinear drift growth. Tamed schemes modify the drift term to control growth without sacrificing explicitness \cite{166}. For instance, a tamed Euler method replaces $f(x)$ by $f(x)/(1+\alpha\Delta t\|f(x)\|)$, ensuring boundedness of increments while preserving strong convergence of order one under appropriate conditions \cite{167}. Implicit schemes, such as backward Euler–Maruyama, provide improved stability properties, particularly for stiff stochastic systems and SPDEs with dissipative operators \cite{168}. Balanced methods interpolate between explicit and implicit updates, adjusting drift and diffusion contributions to simultaneously control stability and accuracy \cite{169}. These approaches have been extended to higher–order stochastic Runge–Kutta methods and to adaptive time–stepping algorithms that refine resolution in regions of high stochastic activity \cite{170}.

\subsection{Long–Time Accuracy and Invariant Measures}
\label{subsec:longtime_invariant}

\begin{figure}[htbp]
    \centering
    \includegraphics[width=0.6\textwidth,keepaspectratio]{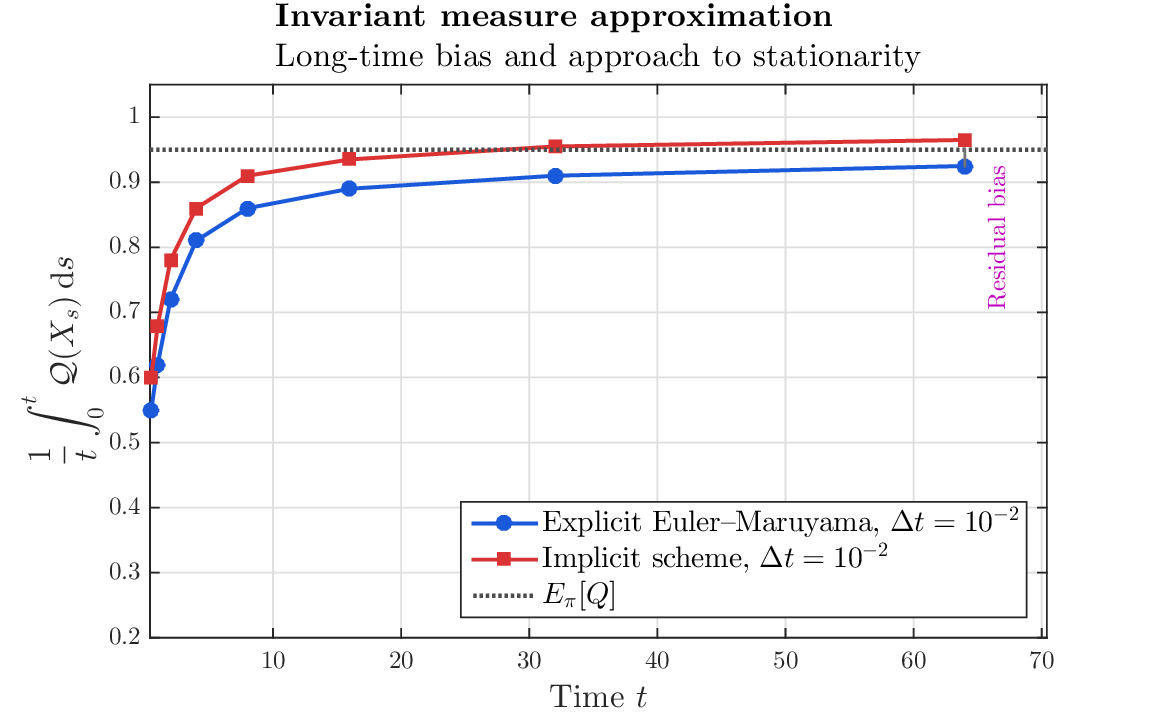}
    \caption{Convergence of time averages to the invariant expectation for an ergodic observable $\mathcal{Q}$. Curves show $\frac{1}{t}\int_0^t \mathcal{Q}(X_s)\,\mathrm{d}s$ for an explicit Euler--Maruyama discretization and an implicit method at the same step size. The dotted line marks $\mathbb{E}_\pi[\mathcal{Q}]$, and the vertical segment at the final time indicates the residual bias. The implicit scheme approaches the stationary value more rapidly and with smaller long–time bias, consistent with theoretical predictions based on generator perturbations and Poisson–equation error analysis.}
    \label{fig:invariant_measure_bias}
\end{figure}

For ergodic SDEs, a central requirement on a numerical method is the faithful approximation of the invariant law over long horizons \cite{171}. If $\pi$ denotes the stationary distribution of the exact process, a method with step size $\Delta t$ induces a Markov chain with stationary law $\pi^{\Delta t}$ \cite{172}. Long–time accuracy entails weak convergence $\pi^{\Delta t}\Rightarrow \pi$ as $\Delta t\to 0$ and a bias in ergodic averages that decays at an optimal rate \cite{173}. The asymptotic bias of an observable $\mathcal{Q}$ admits Poisson–equation representations that expose its dependence on the generator and the local truncation error \cite{174}; this viewpoint underpins modified–equation corrections and Talay–Tubaro–type expansions for invariant–measure error \cite{175}. In practice, drift and diffusion adjustments, symmetrization, and schemes that preserve reversibility or detailed balance reduce long–time bias while retaining acceptable finite–time accuracy \cite{176}. In infinite dimensions, structure–preserving time integrators for SPDEs target invariants and statistical diagnostics dictated by the physics, such as energy spectra \cite{177}, enstrophy \cite{178}, or mass \cite{179}, and are analyzed using spectral gaps or Wasserstein–contractivity of the associated Markov semigroups \cite{180}. Figure~\ref{fig:invariant_measure_bias} illustrates the approach of time averages toward the invariant expectation and highlights the reduced residual bias of an implicit method relative to explicit Euler–Maruyama at a fixed step size.

\subsection{Multilevel Monte Carlo and Rare–Event Simulation}
\label{subsec:mlmc_rare}

Multilevel Monte Carlo reduces work by coupling approximations across a hierarchy in time and, when relevant, space \cite{181}. Let $\varphi_\ell$ denote an approximation of $\varphi(X_T)$ at level $\ell$ with step $\Delta t_\ell$ (and mesh width $h_\ell$ if spatial discretization is present). The telescoping identity:
\[
\mathbb{E}[\varphi_L]=\mathbb{E}[\varphi_0]+\sum_{\ell=1}^{L}\mathbb{E}\big[\varphi_\ell-\varphi_{\ell-1}\big]
\]
is estimated by Monte Carlo averages with $N_\ell$ \emph{coupled} samples \cite{182}. Three rates govern complexity: the weak bias $\mathcal{O}(\Delta t_L^\alpha)$, the variance decay $\operatorname{Var}(\varphi_\ell-\varphi_{\ell-1})=\mathcal{O}(\Delta t_\ell^\beta)$, and the per-sample cost $\mathcal{O}(\Delta t_\ell^{-\gamma})$. With the optimal allocation $N_\ell\propto\sqrt{\operatorname{Var}(\varphi_\ell-\varphi_{\ell-1})/C_\ell}$,
\[
\mathcal{C}_{\text{MLMC}}=
\begin{cases}
\mathcal{O}(\varepsilon^{-2}), & \beta>\gamma,\\[2pt]
\mathcal{O}(\varepsilon^{-2}\log^2\varepsilon), & \beta=\gamma,\\[2pt]
\mathcal{O}\big(\varepsilon^{-2-(\gamma-\beta)/\alpha}\big), & \beta<\gamma,
\end{cases}
\]
which improves on single-level sampling whenever $\beta\ge \gamma$. For Lipschitz SDEs with Euler–Maruyama couplings one typically has $\alpha=1$, $\beta\approx 1$, and $\gamma\approx 1$, leading to $\mathcal{O}(\varepsilon^{-2}\log^2\varepsilon)$ \cite{183}. Antithetic constructions or Milstein-type ideas without Lévy areas often restore $\beta>\gamma$ and the canonical $\mathcal{O}(\varepsilon^{-2})$ \cite{184}. In SPDEs, multilevel couplings combine temporal integrators with spectral or finite-element hierarchies; multi-index Monte Carlo extends the telescoping structure across time–space refinement directions \cite{185}.

\begin{figure}[htbp]
  \centering
  \resizebox{1.05\linewidth}{!}{%
  \begin{tikzpicture}[
      >=Latex,
      font=\sffamily\small,
      line width=1pt,
      node distance=6mm and 8mm, 
      arr/.style   ={->, line width=1pt, draw=Slate},
      darr/.style  ={->, line width=0.9pt, draw=Indigo, dashed},
      card/.style={draw=Slate, rounded corners=3pt, align=center,
                   inner sep=4pt, blur shadow={shadow blur steps=6}},
      level/.style={card, minimum width=36mm, minimum height=9mm,
                    top color=Sky!40, bottom color=white, draw=Slate},
      diff/.style ={card, minimum width=40mm, minimum height=9mm,
                    top color=ElsevierBlue!30, bottom color=white, draw=ElsevierBlue!80!black},
      sum/.style  ={card, minimum width=40mm, minimum height=11mm,
                    top color=Teal!30, bottom color=white, draw=Teal!70!black},
      est/.style  ={card, minimum width=44mm, minimum height=12mm,
                    top color=Indigo!30, bottom color=white, draw=Indigo!80!black},
      rare/.style ={card, minimum width=40mm, minimum height=10mm,
                    top color=Amber!20, bottom color=white, draw=Amber!80!black},
      slim/.style ={font=\sffamily\footnotesize, text=Slate}
    ]

    \node[level] (L0) {$\ell=0$\\\textbf{Coarse path}};
    \node[level,below=10mm of L0] (L1) {$\ell=1$\\\textbf{Coupled paths}};
    \node[level,below=10mm of L1] (L2) {$\ell=2$\\\textbf{Coupled paths}};
    \node[level,below=10mm of L2] (LL) {$\ell=L$\\\textbf{Finest path}};

    \node[diff,right=20mm of L0] (D0) {$\mathbb{E}[\varphi_0]$};
    \node[diff,right=20mm of L1] (D1) {$\mathbb{E}[\varphi_1-\varphi_0]$};
    \node[diff,right=20mm of L2] (D2) {$\mathbb{E}[\varphi_2-\varphi_1]$};
    \node[diff,right=20mm of LL] (DL) {$\mathbb{E}[\varphi_L-\varphi_{L-1}]$};

    \draw[arr] (L0) -- (D0);
    \draw[arr] (L1) -- (D1);
    \draw[arr] (L2) -- (D2);
    \draw[arr] (LL) -- (DL);

    \node[sum, right=16mm of D1] (SUM) {\textbf{Telescoping sum}};
    \draw[arr] (D0) -- (SUM);
    \draw[arr] (D1) -- (SUM);
    \draw[arr] (D2) -- (SUM);
    \draw[arr] (DL) -- (SUM);

    \node[est, right=18mm of SUM] (EST) {\textbf{MLMC estimator}\\$\widehat{\mathbb{E}}[\varphi_L]$};
    \draw[arr] (SUM) -- (EST);

    \node[rare, above=11mm of EST] (IS) {\textbf{Importance sampling}};
    \node[rare, right=10mm of EST] (SPL) {\textbf{Splitting / subset simulation}};
    \node[rare, below=11mm of EST] (AMS) {\textbf{Adaptive multilevel splitting}};

    \draw[darr] (IS)  -- (EST);
    \draw[darr] (SPL) -- (EST);
    \draw[darr] (AMS) -- (EST);

    \node[slim, below=49mm of $(L2)!0.5!(EST)$, anchor=center] (Legend) {%
      \tikz{\node[level, minimum width=8mm, minimum height=3.5mm]{};} \, Hierarchy levels \quad
      \tikz{\node[diff, minimum width=8mm, minimum height=3.5mm]{};} \, Level differences \quad
      \tikz{\node[sum, minimum width=8mm, minimum height=3.5mm]{};} \, Telescoping sum \quad
      \tikz{\node[est, minimum width=8mm, minimum height=3.5mm]{};} \, MLMC estimator \quad
      \tikz{\node[rare, minimum width=8mm, minimum height=3.5mm]{};} \, Rare-event layers
    };

  \end{tikzpicture}%
  }
  \caption{Schematic representation of the { MLMC pipeline} augmented with rare–event variance–reduction strategies. 
On the left, a hierarchy of discretization levels $\ell=0,\dots,L$ (light blue) is constructed, with the coarsest level $\ell=0$ at the top and the finest resolution $\ell=L$ at the bottom. Each level generates either a direct expectation $\mathbb{E}[\varphi_0]$ or a level difference $\mathbb{E}[\varphi_\ell - \varphi_{\ell-1}]$, which collectively form a \emph{telescoping sum}.  This telescoping identity ensures that the computational effort is concentrated where variance is largest, allowing coarser levels to provide inexpensive variance correction to finer ones. 
The sum feeds into the \emph{MLMC estimator}, which delivers an unbiased approximation to $\mathbb{E}[\varphi_L]$ at reduced cost compared to single–level Monte Carlo for the same accuracy. Finally, near the estimator stage, rare–event techniques (including \emph{importance sampling}, \emph{splitting/subset simulation}, and \emph{adaptive multilevel splitting} )can be incorporated to target probabilities of rare events or tail–dependent statistics, further improving estimator efficiency and robustness. 
}
  \label{fig:mlmc_rare_pipeline}
\end{figure}
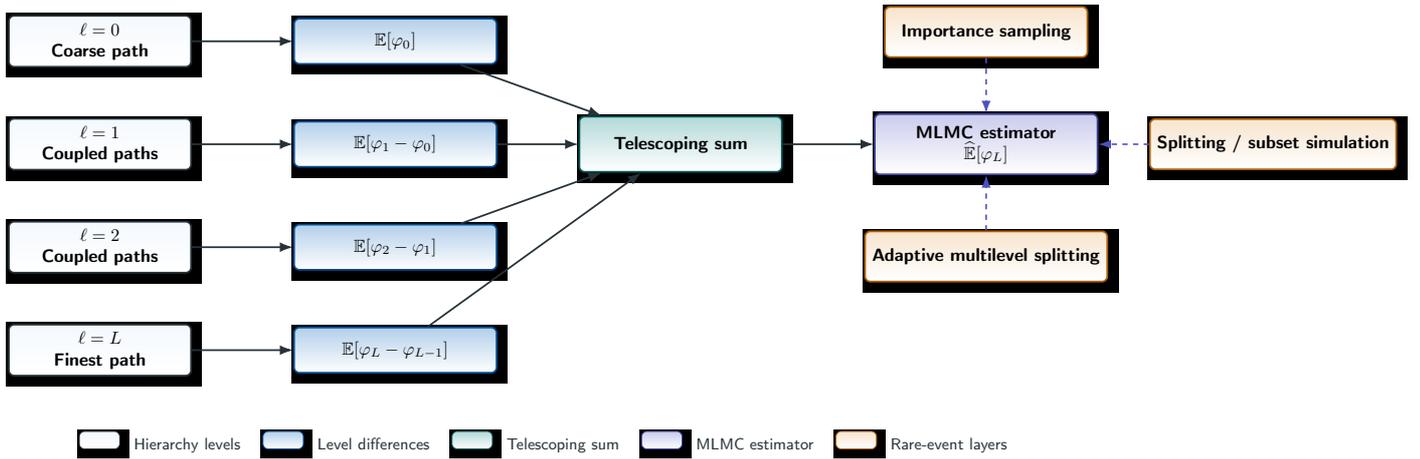

Rare-event probabilities demand variance control beyond stratification by level \cite{186}. Importance sampling with optimally tilted measures derived from large-deviation asymptotics, splitting or subset simulation with adaptively placed intermediate thresholds, and adaptive multilevel splitting or interacting particle systems yield estimators that remain efficient as the probability decreases \cite{187}. These mechanisms integrate naturally with the multilevel decomposition by biasing each difference $\varphi_\ell-\varphi_{\ell-1}$, by coupling biased path pairs, or by embedding control-based surrogates of the optimal change of measure \cite{188}. The pipeline in Fig.~\ref{fig:mlmc_rare_pipeline} summarizes this workflow from hierarchy construction, through telescoping and estimator assembly, to variance-reduction layers for the tail regime.

Overall performance reflects the balance between bias order, coupling variance, long-time statistical fidelity, and variance-reduction design. In practice, structure-preserving or higher-order discretizations are paired with MLMC (or multi-index MLMC) for bulk uncertainty quantification, while importance sampling, splitting, or adaptive multilevel splitting is activated selectively for tail-sensitive observables.

\section{Applications in Biological and Epidemiological Modeling}
\label{sec:applications}

Table~\ref{tab:method_coverage_nocolor} maps method families from Sections~\ref{sec:learning}--\ref{sec:numerics} to concrete uses in biology and epidemiology (2020--2025). We group results into: (A) neural SDE identification and learning; (B) operator learning and physics-informed approaches for PDE/SPDE models; (C) reverse-time generative SDEs; (D) numerical analysis tailored to stochastic models with structural constraints; and (E) inference, surrogates, and phenomena. A cell marked \ok\ indicates at least one relevant published or public preprint; \x\ marks a gap (no result found as of August~2025). The emphasis across domains is consistent: (i) estimating time-varying transmission/reaction rates; (ii) handling partial observability; (iii) preserving positivity and invariants; and (iv) propagating uncertainty for policy-oriented forecasting.

Neural SDEs provide flexible, mechanistically grounded parameterizations for noisy biological dynamics and epidemic time series, learning both drift and diffusion directly from data \cite{189}. Likelihood-based calibration via Euler--Maruyama transitions remains a practical default for discretely sampled outbreaks (e.g., waterborne or contact-driven systems) \cite{200}. Where signal-to-noise and stiffness warrant it, higher-order schemes (Milstein, It\^o--Taylor, SRK) improve statistical efficiency for simulation-and-fit loops \cite{196}. Positive-semidefinite diffusion parameterizations (e.g., Cholesky factors) avoid degeneracies during training in multivariate settings \cite{201}. To tame superlinear drifts that arise from saturating incidence or autocatalysis, monotonicity-preserving or tamed architectures stabilize learning \cite{202}. Continuous-time training through differentiable SDE solvers has been explored in biopharmaceutical and biophysical system identification \cite{203}. Because observation processes are partial and biased in practice, differentiable filters (EnKF/particle) are coupled to neural SDEs for joint state–parameter inference \cite{204}, while martingale/moment methods provide complementary, interpretable checks \cite{205}. Bayesian neural SDEs (e.g., SG-MCMC, variational ensembles) yield posterior trajectories and calibrated intervals for epidemic indicators \cite{206}.

When spatial coupling and nonlocal effects matter, operator learning extends beyond low-dimensional SDEs. Neural operators (Fourier/kernel variants) and PINN-style training have been adapted to SPDE formulations of reaction–diffusion and metapopulation models, improving extrapolation while retaining mechanistic constraints \cite{207,208}. Long-time alignment via ergodic/statistical losses helps match empirical stationary behavior in biochemical and ecological networks \cite{209}. Multirate/multiscale designs address separated time scales (e.g., incubation vs.\ spread; micro vs.\ macro biology) \cite{210}. Hybrid constructions (learned spatial bases with latent temporal SDEs) provide compact surrogates for spatially aggregated epidemics \cite{211}. By contrast, mesh/domain generalization with explicit conservation/positivity guarantees and operator-level Bayesian UQ remain open and underreported (\x).

Score-based diffusion models reinterpret dynamics in reverse time, enabling powerful priors and data augmentation for biological trajectories and epidemic curves \cite{212}. Predictor–corrector samplers support accurate path synthesis and can be adapted for parameter exploration in low-data regimes \cite{213}. Important opportunities remain in controlled reverse-time formulations (links to stochastic optimal control), variance reduction, and reduced-order sampling that respects mechanistic constraints (\x).

Credible inference and forecasting depend on discretizations that respect biology: positivity of compartments, invariance of feasible sets, and stability under stiffness. Recent studies catalogue strong/weak convergence and stability for schemes used in epidemiological SDEs and SPDEs, including linearly implicit and explicit positivity-preserving integrators \cite{195,197,199}. For stochastic reaction–diffusion, unconditionally stable backward–Euler/implicit FD schemes support practical step-size rules \cite{191}. Explicit positivity-preserving methods come with extinction guarantees in subcritical regimes \cite{199}. Underexplored in applied work are tamed/balanced schemes for superlinear drifts, adaptive time stepping, and multilevel Monte Carlo for forward UQ (\x), despite clear relevance to bursty transmission and regime shifts.

Sequential Monte Carlo enables online calibration of latent states and time-varying parameters, handling reporting delays and regime changes \cite{194}. Likelihood-free neural posterior estimation amortizes inference across mechanistic simulators of stochastic epidemics, allowing rapid scenario analysis once trained \cite{190}. Neural surrogates for stochastic SIS-type models reduce evaluation cost in design/control loops \cite{193}. Interacting jump–diffusions and common-noise models illuminate synchronization and chaos-like variability in metapopulations and environmentally forced systems \cite{192}. Stochastic Runge–Kutta integrators round out practical toolkits for vector-borne diseases and multi-scale uncertainty propagation \cite{196}.

Across single-cell dynamics, pharmacological response, vector-borne infections (e.g., dengue), waterborne outbreaks (e.g., cholera), and respiratory epidemics (e.g., COVID-19/monkeypox), four themes recur: (i) \emph{mechanistically informed learning} (neural SDEs, PINNs/operators); (ii) \emph{partial-observation inference} (filters, amortized posteriors); (iii) \emph{structure-preserving numerics} (positivity, invariants, stability); and (iv) \emph{scalable UQ/surrogates} for policy stress-testing. Combining PSD-constrained diffusions, identifiable drift parameterizations, and filter-based training with positivity-preserving or implicit schemes yields the most stable pipelines in practice.


\begin{longtblr}[
  caption = {Recent works (2020--2025) applying the methods surveyed in Sections~\ref{sec:learning}--\ref{sec:numerics} to biology and epidemiology. Legend: \ok\ = at least one study found in the literature; \x\ = no results found to date.},
  label   = {tab:method_coverage_nocolor},
]{%
  colspec = {@{} Q[l,3.6cm] Q[l,2.8cm] Q[c,1.2cm] Q[l,2.0cm] X[3] @{}}, 
  rowhead = 1,
  row{1}  = {font=\bfseries},
  rowsep  = 2pt,
  abovesep= 2pt,
  belowsep= 2pt
}
\toprule
Method (from Secs.~\ref{sec:learning} - \ref{sec:numerics}) & Area & Status &  Refs & Scope / notes \\
\midrule

\SetCell[c=5]{l}\textbf{A. Neural SDE identification and learning} \\
Neural SDEs (learned drift/diffusion) & Neural SDE & \ok & \cite{189} &
scDiffEq learns drift and diffusion from time-resolved single-cell data. \\
MLE with Euler--Maruyama transitions & SDE identification & \ok & \cite{200} &
Applied to cholera epidemic data (discretized likelihood calibration). \\
Higher-order SDE integrators for learning (Milstein / It\^o--Taylor / SRK) & SDE numerics & \ok & \cite{196} &
SRK integrators demonstrated on dengue SDEs (sampling / simulation focus). \\
PSD diffusion via Cholesky factorization & SDE identification & \ok & \cite{201} &
Positive-semidefinite diffusion parameterization in biological stochastic systems. \\
Tamed / monotonicity-preserving architectures & Robust learning & \ok & \cite{202} &
Network-informed generative modeling for epidemic prediction / reconstruction. \\
Continuous-time training through SDE solvers & System ID (biopharm.) & \ok & \cite{203} &
Continuous-time echo state networks for pharmacological dynamics. \\
Partially observed training with differentiable filters (EnKF / particles) & State--parameter inference & \ok & \cite{204} &
Epidemiologically informed particle filtering (monkeypox, 2022). \\
Moment-based (martingale) estimation & Stochastic processes & \ok & \cite{205} &
Population CTMCs; bounds on mean first passage times (moment equations). \\
Bayesian neural SDEs (SG-MCMC / ensemble VI) & Bayesian inference & \ok & \cite{206} &
Posterior inference for COVID-19 dynamics in Germany. \\
Bootstrap UQ over trajectories / noise seeds & Uncertainty quantification & \x & -- &
No results in the literature. \\

\SetCell[c=5]{l}\textbf{B. Operator learning and PINNs for SPDE / PDE models} \\
Neural operators (Fourier / kernel) for SPDEs & Operator learning (SPDE) & \ok & \cite{207} &
Physics-informed neural treatments for fractional epidemiological models. \\
Physics-Informed Neural Operators / PINNs & Operator learning & \ok & \cite{208} &
PINN-style training for compartmental epidemiological dynamics. \\
Ergodic / statistical losses for long-time alignment & Operator learning & \ok & \cite{209} &
Ergodicity constraints guide inference in latent biological reaction dynamics. \\
Mesh / domain generalization and conservation constraints & Operator learning (SPDE) & \x & -- &
No results in the literature. \\
Multirate / multiscale operator architectures & Operator learning & \ok & \cite{210} &
Multi-level stochastic modeling framework for computational epidemiology. \\
Operator-level UQ (Bayesian priors, ensembles, linearization) & Operator learning (SPDE) & \x & -- &
No results in the literature. \\
Hybrid latent SDE with learned spatial bases & Hybrid SDE--SPDE & \ok & \cite{211} &
Low-rank latent temporal dynamics for COVID-19 (hybrid representation). \\

\SetCell[c=5]{l}\textbf{C. Generative SDEs and reverse-time modeling} \\
Reverse-time SDE and probability-flow ODE & Generative SDEs & \ok & \cite{212} &
Reverse-time inference for biological dynamics (score-based viewpoint). \\
Predictor--corrector samplers for reverse SDE & Generative SDEs & \ok & \cite{213} &
Parameter estimation strategies for nonlinear biological models. \\
Adaptive stepping / variance reduction / reduced-order sampling & Generative SDEs & \x & -- &
No results in the literature. \\
Controlled reverse-time (SDE $\leftrightarrow$ stochastic control) & Generative SDEs & \x & -- &
No results in the literature. \\

\SetCell[c=5]{l}\textbf{D. Numerical analysis for SDEs / SPDEs} \\
Strong / weak convergence and stability frameworks & Numerics (SDE/SPDE) & \ok & \cite{195,197,199} &
Monograph baseline; linearly implicit SIS; positivity-preserving explicit schemes. \\
Implicit schemes for SPDEs (backward--Euler + implicit FD) & SPDE numerics & \ok & \cite{191} &
Unconditionally stable schemes for stochastic reaction--diffusion in epidemics. \\
Linearly implicit schemes for SIS & SDE numerics & \ok & \cite{197} &
Stability regions; accuracy trade-offs for stiff nonlinear incidence. \\
Tamed schemes for superlinear drift & SDE numerics & \x & -- &
No results in the literature. \\
Balanced schemes & SDE numerics & \x & -- &
No results in the literature. \\
Adaptive time stepping for SDE / SPDE & SDE/SPDE numerics & \x & -- &
No results in the literature. \\
Long-time accuracy / invariant-measure fidelity & SDE/SPDE numerics & \ok & \cite{195} &
Survey-level treatment of invariant measures / long-horizon bias. \\
Structure / constraint preserving (positivity / invariance) & SDE/SPDE numerics & \ok & \cite{199} &
Explicit positivity-preserving scheme with extinction guarantees. \\
Multilevel Monte Carlo / multi-index MLMC & Uncertainty quantification & \x & -- &
No results in the literature. \\

\SetCell[c=5]{l}\textbf{E. Inference, surrogates, and phenomena} \\
Sequential Monte Carlo for online calibration & Inference & \ok & \cite{194} &
Joint tracking of latent states and time-varying parameters (regime shifts, bias). \\
Likelihood-free NPE for stochastic epidemics & Inference & \ok & \cite{190} &
Amortized posterior estimation for mechanistic epidemic models. \\
ANN surrogates for stochastic SIS & Surrogates & \ok & \cite{193} &
LM-trained emulator for nonlinear SIS with stochasticity. \\
SPDE chaotic variability with jumps / common noise & Modeling phenomena & \ok & \cite{192} &
Interacting jump--diffusion (McKean--Vlasov with common noise; chaos-like signatures). \\
Stochastic Runge--Kutta for epidemic SDEs & SDE numerics & \ok & \cite{196} &
SRK integrators demonstrated for dengue uncertainty. \\
\bottomrule
\end{longtblr}

\section{Outlook and Open Problems}
\label{sec:outlook}

The integration of data-driven system identification, continuous-time generative modeling, and structure-aware numerics is transforming stochastic modeling in biology and epidemiology. Table~\ref{tab:open_problem_matrix_refined} outlines the most immediate opportunities, detailing concrete goals, candidate methods, and measurable evaluation criteria.

For neural SDEs under partial observation and state-dependent diffusion, we need verifiable identifiability and finite-sample rates that prevent drift–diffusion confounding and tolerate irregular sampling. For operator learning, mesh-independent approximation and stability guarantees should extend to random forcings, with conditions ensuring strong solutions, pathwise stability, geometric ergodicity, and (when Bayesian) posterior contraction. \emph{Why it matters:} principled recovery of dynamics and uncertainty rather than ad hoc fitting (see Table~\ref{tab:open_problem_matrix_refined}, rows “Identifiability’’ and “Consistency’’).

Training objectives should be \emph{discretization-aware}: optimize under the same integrator used at inference to avoid train–test drift. Memory-stable adjoints and AD-compatible strong/stable schemes (tamed/balanced Milstein, linearly implicit IMEX, positivity-preserving updates) are needed for stiff SDE/SPDE regimes; adaptive estimators should separate discretization from statistical error \emph{inside} training loops (Table~\ref{tab:open_problem_matrix_refined}, “Solver-aware’’ and “Adaptive error control’’).

Move beyond pointwise error bars to calibrated, shift-robust sequential forecasts: decompose aleatoric vs.\ epistemic uncertainty, use dependence-aware conformal layers, and report distributional scores (CRPS/energy), PIT diagnostics, and long-horizon bias of invariant statistics. Tail risk remains a gap: couple multi-index MLMC with rare-event biasing to obtain predictable variance–cost behavior.

\begin{table*}[htbp]
\centering
\caption{Open problems distilled into goals, candidate methods, and reportable metrics for learned stochastic dynamics and structure-aware numerics.}
\label{tab:open_problem_matrix_refined}
\renewcommand{\arraystretch}{1.12}
\begin{tabularx}{\textwidth}{@{}L{3.7cm} X X L{3.0cm}@{}}
\toprule
\textbf{Theme} & \textbf{Goal} & \textbf{Candidate methods} & \textbf{Metrics}\\
\midrule
Identifiability (neural SDE) &
Disentangle drift vs.\ diffusion under partial observation and irregular sampling. &
Structured diffusion (factorized/low-rank), short-increment penalties; stationary-moment constraints; uniform empirical-process bounds; confounding tests. &
Drift/diffusion MSE; misattribution rate; CI coverage. \\
\addlinespace[0.8mm]
Consistency \& well-posedness &
Guarantee strong solutions, stability, ergodicity, and (Bayesian) posterior contraction. &
One-sided Lipschitz \& linear-growth priors; Lyapunov drift/minorization; Wasserstein contractivity; sieve priors and small-ball conditions. &
Stability constants; mixing rate; ergodic bias; contraction rate. \\
\addlinespace[0.8mm]
Solver-aware training &
Remove train–inference discretization mismatch. &
Discretization-aligned losses; reversible/checkpointed adjoints; step-size–aware objectives; gradient-stability control for stiffness. &
Generalization gap (train vs.\ inference); gradient growth; solver-induced bias. \\
\addlinespace[0.8mm]
Structure-preserving SPDE integrators &
Maintain invariants (mass/positivity) and spectra over long horizons. &
Linearly implicit IMEX; positivity-preserving projections; energy/enstrophy-aware updates; invariant-region enforcement. &
Invariant error; spectral discrepancy; long-time bias. \\
\addlinespace[0.8mm]
Adaptive error control (in-loop) &
Separate discretization vs.\ statistical error and adapt time/mesh during training. &
Dual-weighted residuals; MLMC/MIMC gradient estimators; online mesh/time refinement. &
Work–RMSE tradeoff; estimator variance; speedup at fixed error. \\
\addlinespace[0.8mm]
Tail risk \& rare events &
Estimate exceedance probabilities with predictable variance–cost. &
Schr\"odinger-bridge or Girsanov tilting; adaptive splitting/AMS; multi-index MLMC couplings. &
Relative error at fixed cost; variance-decay exponent; ESS. \\
\addlinespace[0.8mm]
Benchmarking \& reporting &
Make results comparable and auditable across methods/costs. &
Fixed data splits; structural baselines; standardized metric suite; compute/energy reporting. &
Reproducible leaderboard; compute-normalized scores; energy/carbon budget. \\
\bottomrule
\end{tabularx}

\vspace{4pt}
\footnotesize
\textit{Abbreviations.} MSE: mean-squared error. CI: confidence interval. IMEX: implicit–explicit splitting. MLMC/MIMC: multi(level/index) Monte Carlo. CRPS: continuous ranked probability score. PIT: probability integral transform. AMS: adaptive multilevel splitting. ESS: effective sample size.
\end{table*}


\section{Conclusion}
\label{sec:conclusion}

This review set out to give researchers a compact, navigable map of recent (2020--2025) advances in stochastic modeling for biology and epidemiology, with two goals: (i) to locate and clarify \emph{novelty} across theory, inference, learning, and numerics; and (ii) to surface \emph{open problems} and underused methods with clear relevance to biological and epidemiological practice. We organized the field along a pipeline from Itô SDE/SPDE formulations to neural SDE identification and operator learning, reverse-time generative SDEs, and structure-preserving integrators, so that readers can see where each contribution fits and how choices in representation, solver, and inference jointly determine scientific reliability.

Concretely, we (a) distilled an identification toolkit for drift–diffusion learning under partial observation (with PSD-constrained diffusion, monotonicity/taming, and differentiable filtering), (b) summarized operator-learning surrogates for spatial dynamics with weak-physics and ergodic/statistical alignment, and (c) reviewed solver-aware computation that maintains feasibility, stability, and long-time statistics while enabling scalable uncertainty quantification. We also charted gaps that invite immediate work in biology/epidemiology: identifiability under irregular sampling and state-dependent noise; operator learning under random forcing with well-posedness/ergodicity guarantees; discretization-aware training with memory-stable adjoints; positivity- and invariant-preserving schemes in routine pipelines; multilevel/multi-index Monte Carlo for tail-sensitive policy metrics; controlled reverse-time SDEs for assimilation and counterfactuals; and operator-level Bayesian UQ with mesh/domain generalization.

By mapping what is new, what is missing, and where methods align with biological and epidemiological objectives (time-varying transmission, spatial spread, regime shifts, and decision-oriented UQ), this review aims to shorten the path from method selection to credible results, reduce reinvention, and focus effort on the most impactful open directions.
\section*{Declaration of competing interest}
There is no financial or non-financial assistance provided by a third party for the reported work.
\section*{Data availability}
No data was used for the research described in the article. 


\end{document}